\documentclass{amsart}

\usepackage{graphics}
\usepackage[dvips]{graphicx}
\usepackage[latin1]{inputenc}

\newtheorem{thm}{Theorem}[section]

\theoremstyle{definition}
\newtheorem{dfn}[thm]{Definition}
\theoremstyle{remark}
\newtheorem{rk}[thm]{Remark}
\numberwithin{equation}{section}
\newcommand{\R}{\mathbb R}

\newcommand{\eps}{\varepsilon}

\begin{document}

\author{Vassily Olegovich Manturov}

\address{V.O.Manturov, Chelyabinsk State University}\footnote{This research is supported by the Russian Science Foundation, Project No. 16-11-10291}%
\email{vomanturov@yandex.ru}%

\title{A Note on a Map from Knots to $2$-Knots}

\maketitle

\section{Introduction. Basic Notions}

The aim of the present paper is to construct a well defined map $\alpha$ from classical knots to (some generalization of) $2$-knots, which naturally works for any classical knots in codimension $2$. This map works well for braids and tangles; besides, it is natural with respect to the boundary operation. In particular, if a classical knot $K$ is {\em slice} (i.e., can be spanned by a $2$-dimensional object in four-dimensional half-space bounded by the $3$-space where $K$ lives) then the corresponding object $\alpha(K)$ can also be bounded by a some trivial $2$-dimensional object. This leads to sliceness obstructions for classical knot.

This work naturally is closely connected with \cite{Great}, where various invariants of {\em dynamical systems} valued in {\em groups} were constructed. The main example from \cite{Great}, see also \cite{MN}, is the construction of an invariant of classical braids valued in a certain group $G_{n}^{3}$. Braids are considered as motions of points on the plane; the generators of this group correspond to those moments where some three points are collinear (form a horizontal quadrisecant); in some sense the abstract $2$-knots (diagrams modulo moves) constructed here represent the ``knot'' counterpart, whence group elements (words modulo relations) are the braid counterpart.

\begin{dfn}
A knot $K\subset \R^{3}$ is called {\em slice} if there is a $2$-manifold $M^{2}$ with boundary properly embedded in $\R^{4}_{+}$ such that the image of $\partial M^{2}$ coincides with $K\subset\R^{3}=\partial \R^{4}_{+}$.
\end{dfn}

Now, we pass to {\em abstract surface knot diagrams}, see,e.g.,\cite{KK}.

\begin{dfn}
Let $F$ be a $2$-surface; a {\em double decker set} is a $1$-complex $D\subset F$ together with a pasting rule such that:
\begin{enumerate}
\item the pasting is an equivalence relation for $D$ which is continuous with respect to the topology of $F$;

\item each equivalence class consists of one, two or three points from $D$; these sets are denoted by $D_{1},D_{2},D_{3}$, respectively, points from $D_{2}$ are called {\em double points}.

\item $Card(D\backslash D_{2})$ is finite; $D_{1}\sqcup D_{3}\subset {\bar D_{2}}$;

\item elements of each equivalence class are ordered; the ordering is continuous with respect to the topology of $F$.
For each two pasted points $d_{1},d_{2}$ such that $d_{1}>d_{2}$ we say that $d_{1}$ {\em lies over} $d_{2}$;

\item let $(d_{1},d_{2})$ be a {\em double point} of $D$, say, $d_{1}<d_{2}$; then in the neighbourhoods of $d_{1}$ and $d_{2}$ the double decker set represents two open intervals; these intervals consist of double points; points from one interval are identified with points from the other interval.

\end{enumerate}

\end{dfn}

\begin{dfn}
For a surface $F$ with a double decker set we naturally get a $2$-complex $F^{\sim}$ obtained by identifying all equivalent points.

Such a $2$-complex endowed with additional over/under structure is called an {\em abstract $2$-surface diagram}.
\end{dfn}

\begin{rk}
From the definition above it follows that $D_{1}$ consists of {\em cusps}; a neighbourhood of a cusp topologically looks like a circle $|z|<1$ where $z$ is identified with $-z$ with $z=0$ being the cusp point; similarly, triple points topologically look like the intersection of three hyperplanes in $\R^{3}$.

As for double points, they naturally form {\em double lines} which are $1$-manifolds possibly having boundary at cusps and triple points (as $K_{1}\sqcup K_{3}\subset {\bar K_{2}}$).
\end{rk}

An abstract knot diagram is not assumed to be embedded anywhere. Nevertheless, a local neighbourhood of any point of $F$ can be embedded in $\R^{3}$. Thus, it is natural to depict local parts of an abstract $2$-surface in a way similar to projections of $2$-knots in $3$-space. They naturally appear as an  general position orthogonal projection of a $2$-knot in $\R^{4}$ to some $3$-space; points having more than one preimage are identified; the partial order relation is defined with respect to the projection coordinate.

Usually, we deal with {\em $2$-surface knot diagrams} when $F$ is connected and {\em $2$-surface link diagrams} when $F$ consists of connected components $F=F_{1}\sqcup F_{2}\sqcup \cdots \sqcup F_{n}$; here we say that $n$ is the number of components.

Roseman moves \cite{Ros} are initially defined for $2$-knots in $3$-space.

 However, they are {\em local} and can be applied to abstract knot diagrams.

\begin{figure}[h]
\begin{center}
\includegraphics[width=90mm]{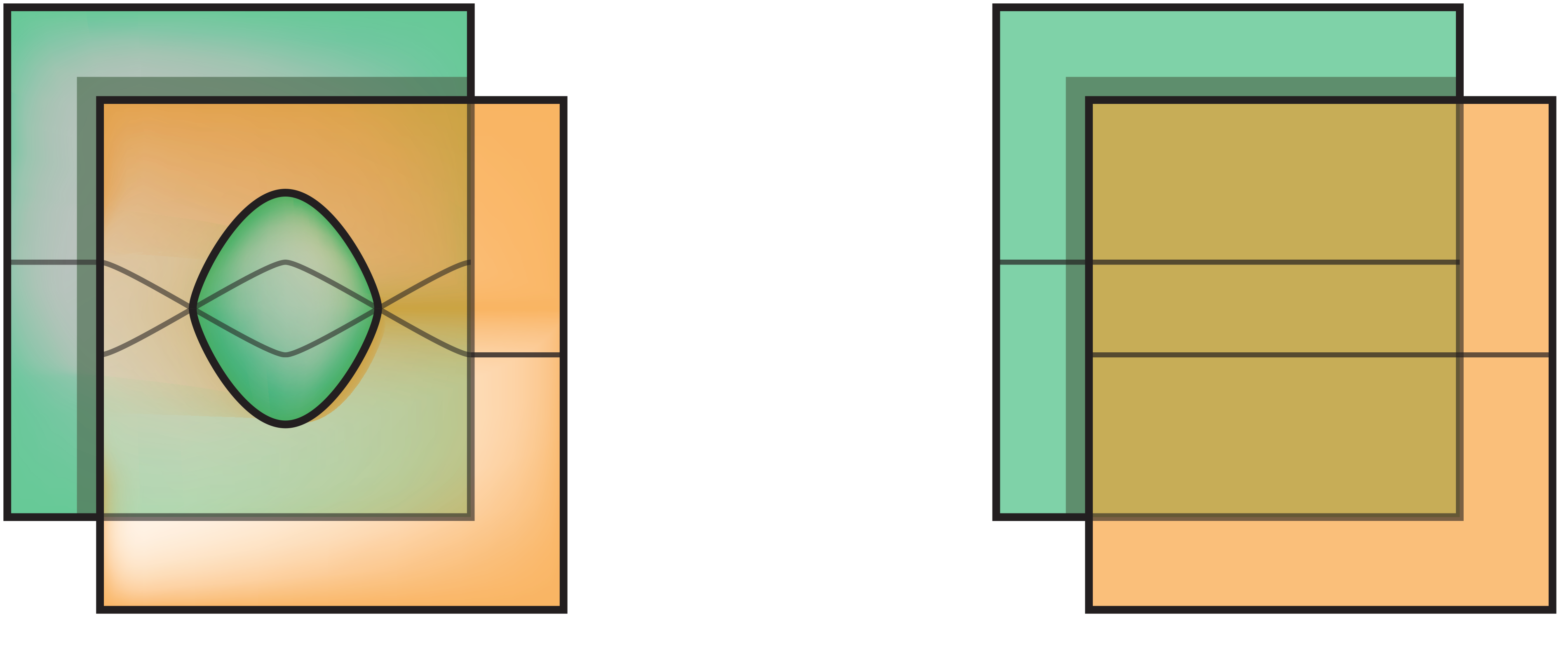}
\caption{The first Roseman move $\mathcal{R}_1$: elliptic move of type $\Omega_2$.} \label{ris:R1}
\end{center}
\end{figure}

\begin{figure}
\begin{center}
$\includegraphics[width=90mm]{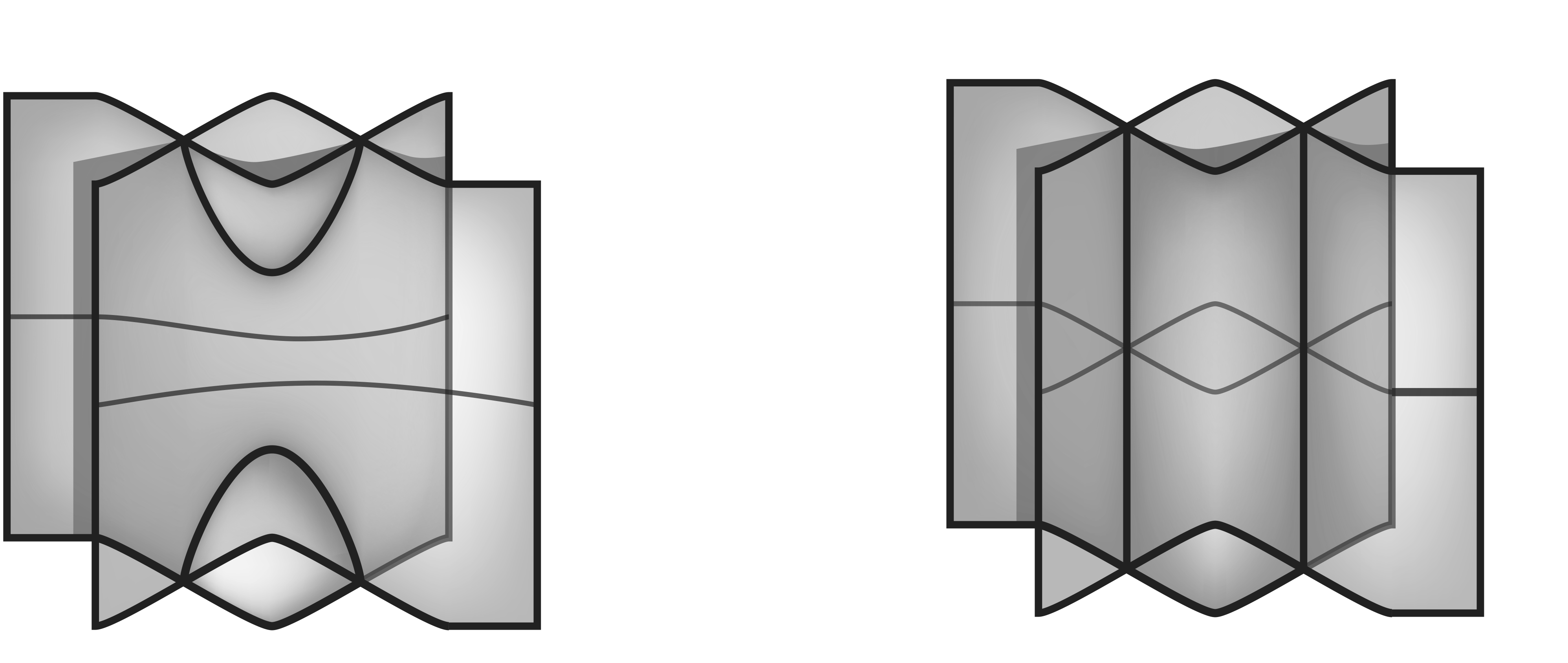}$
\caption{The second Roseman move $\mathcal{R}_2$: hyperbolic type $\Omega_2$ move.} \label{ris:R2}
\end{center}
\end{figure}

\begin{figure}[h]
\begin{center}
$\includegraphics[width=90mm]{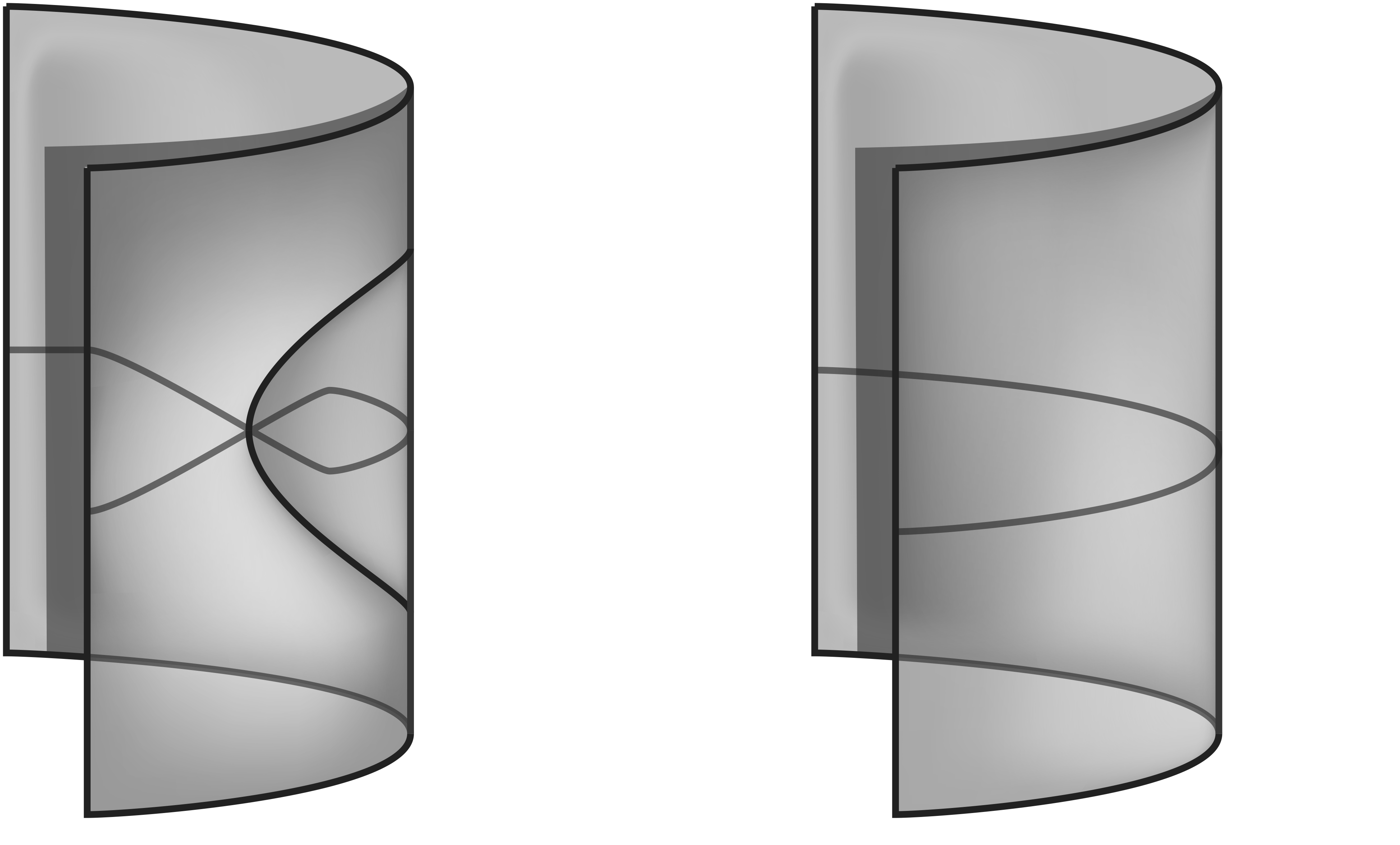}$
\caption{The third Roseman move $\mathcal{R}_3$: elliptic $\Omega_1$ move.} \label{ris:R3}
\end{center}
\end{figure}

\begin{figure}[h]
\begin{center}
$\includegraphics[width=90mm]{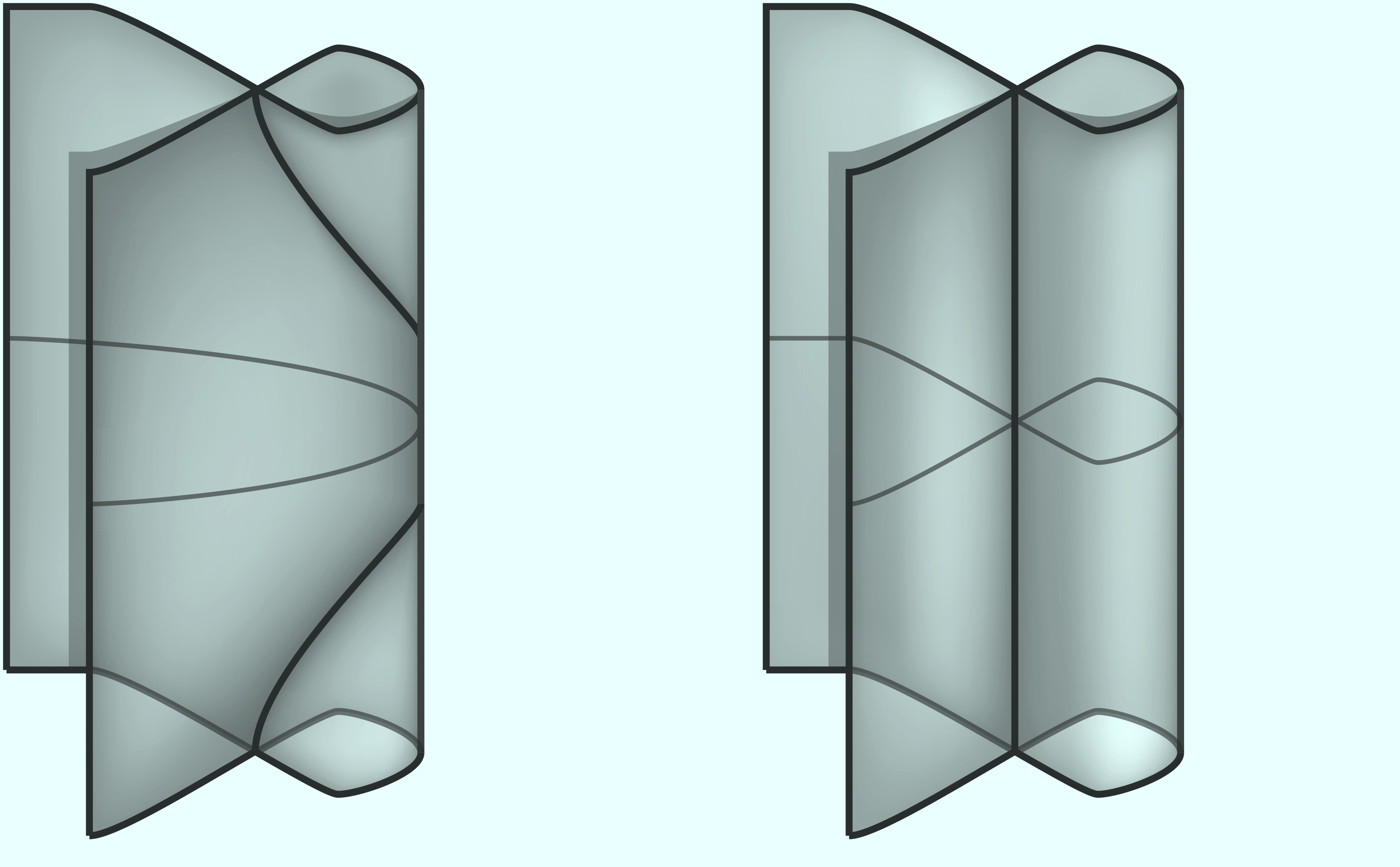}$
\caption{The fourth Roseman move $\mathcal{R}_4$: hyperbolic $\Omega_1$ move.} \label{ris:R4}
\end{center}
\end{figure}

\begin{figure}[h]
\begin{center}
$\includegraphics[width=90mm]{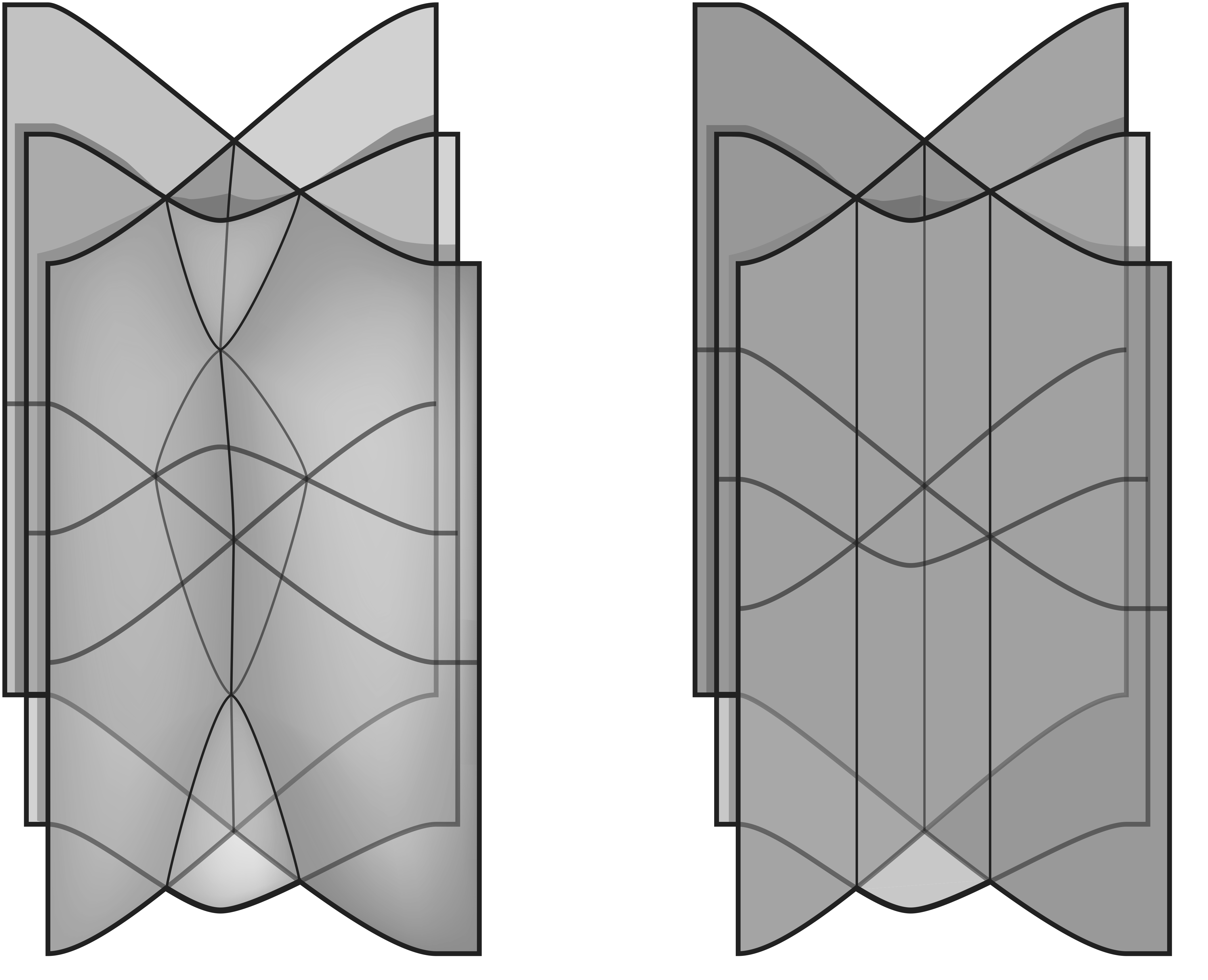}$
\caption{The fifth Roseman move $\mathcal{R}_5$:  $\Omega_3$ move.} \label{ris:R5}
\end{center}
\end{figure}

\begin{figure}[h]
\begin{center}
$\includegraphics[width=90mm]{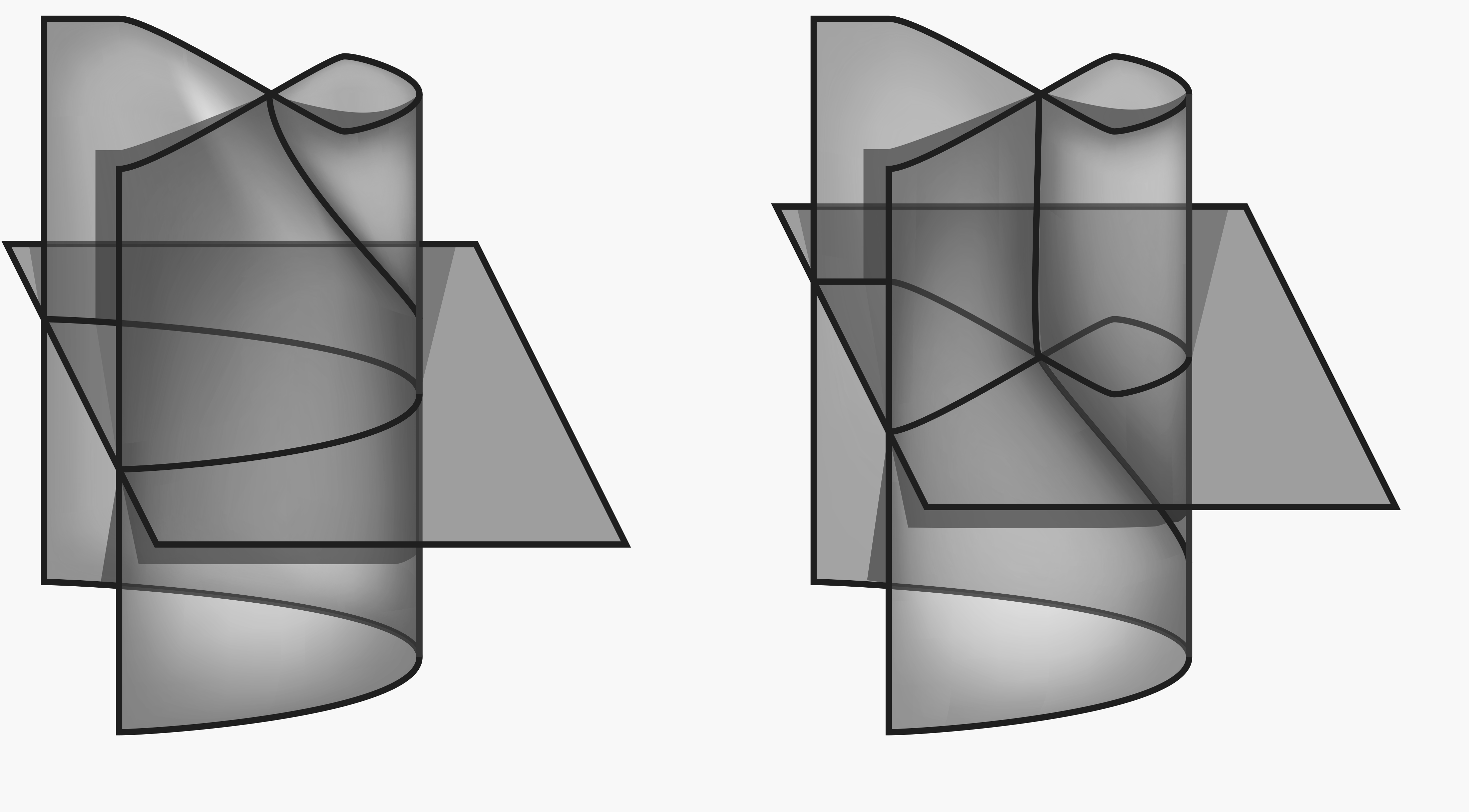}$
\caption{The sixth Roseman move $\mathcal{R}_6$: crossing the branching point.} \label{ris:R6}
\end{center}
\end{figure}

\begin{figure}[h]
\begin{center}
$\includegraphics[width=90mm]{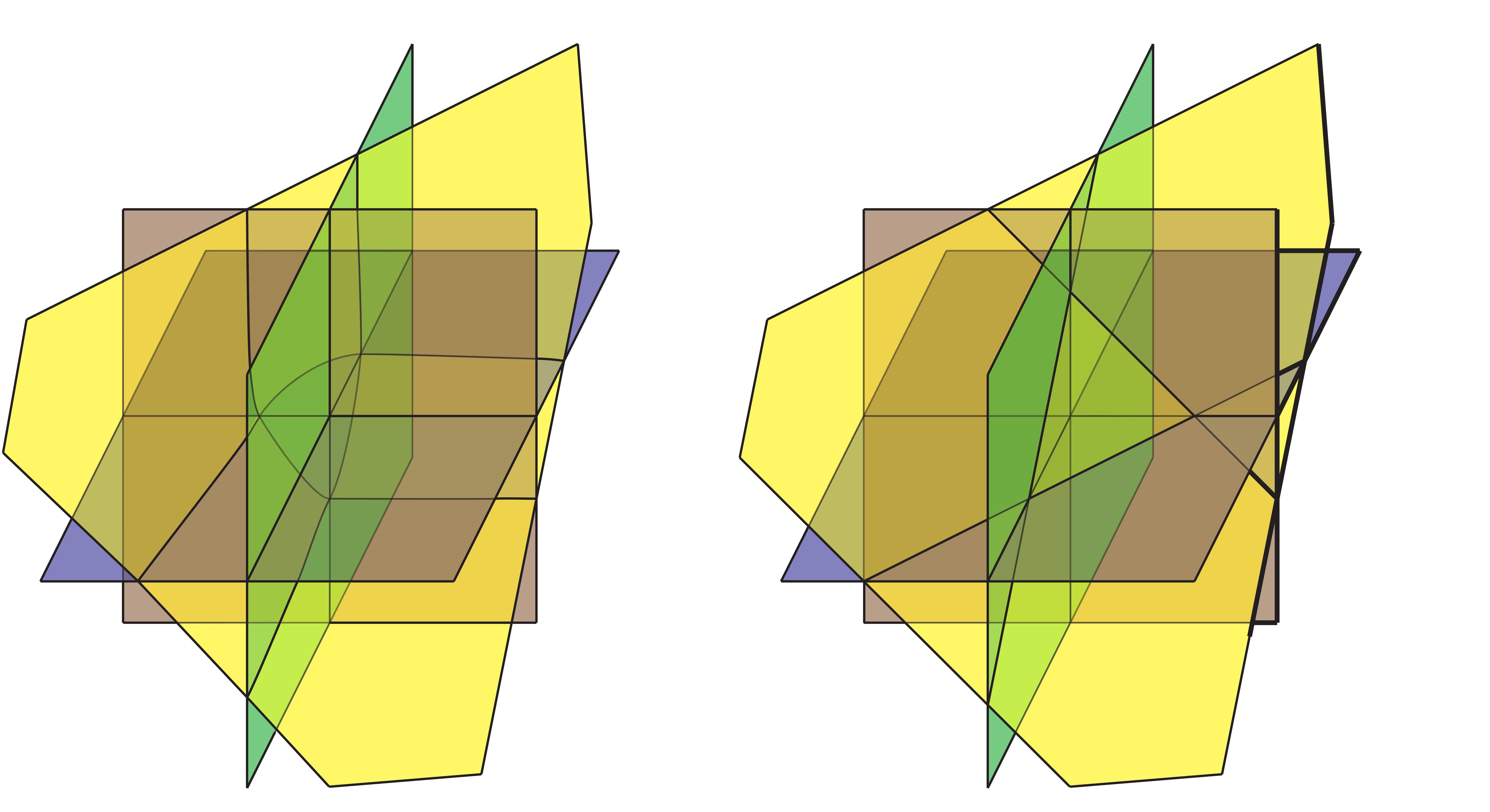}$
\caption{The seventh Roseman move $\mathcal{R}_7$: the tetrahedral move.} \label{ris:R7}
\end{center}
\end{figure}

\begin{thm}[Roseman \cite{Ros}]
Two $2$-knot (resp., $2$-link) diagrams $D,D'$ represent isotopic $2$-knots (resp., $2$-links) if and only if one can get from $D$ to $D'$ by a finite sequence of {\rm Roseman moves}, shown in Figs. \ref{ris:R1}--\ref{ris:R7}, and isotopies in $\mathbb{R}^3$.
\end{thm}

Later on, we make no difference between combinatorially equivalent (abstract) knot diagrams.

Instead of giving a definition of a $2$-knot diagram, we shall pass to a more general definition (for abstract knots).

\begin{dfn}
An {\em abstract surface-knot} is an equivalence class of abstract surface diagrams modulo Roseman moves. Namely, two diagrams $D,D'$ are called {\em equivalent} if there is a sequence of diagrams $D=D_{1}\to D_{2}\to \dots \to D_{n}=D'$, such that each two diagrams $D_{i},D_{i+1}$ are related by a Roseman move.
\end{dfn}

\section{The case of $1$-Knots}

Let $S^{1}_{\phi}$ be the angular unit sphere whose points are identified to unit vectors $e^{i\phi}$. Sometimes we also use the notation $-\phi$ for $\phi+\pi$.

Let $K$ be a knot in $\R^{3}$. Later on, we shall require some {\em general position for $K$}. We shall be especially interested in the height function for $K$; whenever mentioning ``minimum'' or ``maximum'' we assume extrema with respect to this function.
Let $S(K)= T^{2}=S^{1}\times S^{1}_{\phi}$ be the torus where the first coordinate corresponds to the knot $K$ and the second coordinate is $\phi$ for $S^{1}=S^{1}_{\phi}$.

Now, we make an {\em abstract $2$-knot diagram} $D(K)$ of $S(K)$, as follows. We define the equivalence for pairs of points  $(k,\phi),(k',\phi),k\in K,\phi\in S^{1}$ that $k$ and $k'$ are points on $K$ lying on the same horizontal plane, and $\phi$ is such that the pointing vector from $k$ to $k$ is collinear with $\phi$ (or $-\phi$); hence, if $(k,\phi)\sim (k',\phi)$ then $(k,-\phi)\sim (k',-\phi)$ and vice versa; moreover, if $(k,\phi)$ {\em lies over} $(k',\phi)$ then $(k,-\phi)$ lies under $(k',-\phi)$.
We are going to identify the equivalent points;

 \begin{dfn} We say that $K\in \R^{3}$ is in
 {\em general position condition} if the above  identification leads to a double decker set. Namely, we require that
\begin{enumerate}

\item No four points $(k_{1},\phi),(k_{2},\phi),(k_{3},\phi),(k_{4},\phi)$ are equivalent. In other words, this means that $K$ has no horizontal quadrisecant.

\item The number of triples $(k_{1},\phi)~(k_{2},\phi)~(k_{3},\phi)$ is finite; this means that $K$ has finitely many horizontal quadrisecants.

\item There is a finite number of {\em cusps};

\item All other double points $X_{1}=(k_{1},\phi)\sim (k_{2},\phi)=X_{2}$ are {\em regular} in the following sense. Topologically, the set of double points in the neighbourhood of $X_{1}$ is homeomorphic to an open interval as well as the set of double points in the neighbourhood of $X_{2}$, and these two intervals are identified.

\end{enumerate}
\end{dfn}

\begin{rk}
Note that regular points can correspond to {\em local extrema}. For instance, assume $K$ locally looks like the graph of the parabola $z=x^{2}, y=const$ and another branch of $K$ looks just like a vertical line $x=y=1$. Let $X=((0,0,0),\frac{\pi}{4})$ be origin with the vector pointing to the point $(1,1)$. Then the neighbourhood of $X$ consists of two parts, the one with $(\eps,0,\eps^{2})$ for $\eps>0$ and the one for $\eps<0$; in both cases the vector pointing to $(0,0,\eps^{2})$ is close to $\frac{\pi}{4}$. Then $D_{2}$ will locally consist of two double lines, one of which corresponds to the neighbourhood of $\frac{\pi}{4}$ and the other one corresponds to the neighbourhood of $-\frac{\pi}{4}$.
\end{rk}

By $\alpha(K)$ we denote the abstract $2$-knot represented by $D(K)$.

Now, we are ready to formulate the main theorem.
\begin{thm}
If knots $K$ and $K'$ are isotopic in $\R^{3}$ then the abstract surface $2$-knots $\alpha(K)$ and $\alpha(K')$ are equivalent.
\end{thm}

The proof consists of case-by-case consideration of all situations where $K$ fails to be generic.
As an example, we say that if we pass through a horizontal quadrisecant passing through some four points $a,b,c,d\in K$ with the same $z$-coordinate, then yields two ``antipodal'' seventh Roseman moves, see Fig.\ref{ris:R7}. Indeed, denoting the angles of the quadrisecant by $\phi$ and $-\phi$, we can see four sheets of the surface $D(K)$ in the neighbourhoods of $(a,\phi),(b,\phi),(c,\phi),(d,\phi)$. These four sheets have four double points, and the effect of the passing through the quadrisecant is the inversion of the tetrahedron.

The same happens for the direction $-\phi$.

Now, the same definition and the same proof generalizes for many other objects: braids, tangles, and even $2$-knots and knots in higher dimensions. For tangles (or knots in $\R^{2}\times [0,\infty)$) we assume $K$ to be a properly embedded (into $[0,1]$ or $[0,\infty)$) collection of closed intervals; the resulting object $\alpha(K)$ will be an equvalence class of corresponding $2$-surface diagrams with boundary; braids are a special case of tangles.

\section{Dimension $2$}

More importantly, the same definition works in higher dimensions.

Let $K$ be a $2$-knot in $3$-space. Here, we are interested in planes $z=const,t=const$. As before, we let $S(K)= K\times S^{1}_{\phi}$, where the first coordinate corresponds to the knot $K$ and the second coordinate is $\phi$ for $S^{1}=S^{1}_{\phi}$.

For each point $x\in K$, we consider the plane $z=const, t=const$ passing through this point. In general position, for each plane $P=\{z=c_{1},t=c_{2}\}$, the intersection $K\cap P$ consists of finitely many points.

We are not going to give the formal definition of the abstract surface $3$-diagram and an abstract surface $3$-knot (link). They are quite similar to the $2$-dimensional case. However, one should point out that the set of Roseman moves in higher dimension which originates from codimension $1$ singularities for $3$-knot projections to $4$-space, is finite and all such moves are local (The same works in any dimension). As in the $2$-dimensional case, {\em abstract $3$-knot} is the equivalence class of diagrams modulo moves.

As before, we identify those points $(k,\alpha)$ and $(k',\alpha)$ if $k$ and $k'$ have the same third and fourth coordinate.

We can think of this map as follows. Having a $2$-knot $K$, we can consider its slices $K_{c}=K\cap \{t=c\}$; for generic $c$, these slices are just $1$-links in $3$-space; we can think of $K$ as ``glued'' of all $K_{c}$ for distinct $c$. Thus, $\alpha(K)$ is glued in the same way from distinct $\alpha(K_{c})$.

\begin{thm}
If $2$-knots $K$ and $K'$ are isotopic in $\R^{4}$ then the abstract $3$-knots $\alpha(K)$ and $\alpha(K')$ are equivalent.
\end{thm}

The main observation is that for each neighbourhood $(\phi-\eps,\phi+\eps)$ the picture of $\alpha(K)$ can be represented as a collection of sheets in the Euclidean space. Thus, whenever we perform some isotopy of $K$ and restrict, every combinatorial diagram can be drawn in the Eucludean space and is subject to local moves in these Eucliedan space (for dimension $2$ these are exactly Roseman moves).

\begin{thm}
If $K$ is slice then $\alpha(K)$ is slice.
\end{thm}

The above generalization for the case of $2$-knots works as well for any $n$-manifold $M^{n}$ embedded in some fixed space of $\R^{2}\times N^{n}$ for some manifold $N^{n}$; the only thing we need here is the existence of the ``first'' two coordinates $x,y$ so that we can take ``horizontal'' slices $\R^{2}\times \{*\}$ which intersect $M^{n}\subset \R^{2}\times N^{n}$ at finitely many points.

\end{document}